\documentclass{elsarticle}

\usepackage{amsmath}
\usepackage{bm}
\usepackage{tikz}
\usepackage{pgfplots}

\usepackage{hyperref}

\newtheorem{thm}{Theorem}

\newproof{pf}{Proof}

\newcommand{\const}{\mathop{\rm const}\nolimits}

\journal{arXiv.org}

\begin{document}

\begin{frontmatter}

\title{Numerical solving unsteady space-fractional problems
with the square root of an elliptic operator}

\author[nsi,univ]{Petr N. Vabishchevich}
\ead{vabishchevich@gmail.com}

\address[nsi]{Nuclear Safety Institute, Russian Academy of Sciences, Moscow, Russia}
\address[univ]{North-Eastern Federal University,  Yakutsk, Russia}

\begin{abstract}
An unsteady problem is considered for a space-fractional equation in a bounded domain. 
A first-order evolutionary equation involves the square root
of an elliptic operator of second order. Finite element approximation in
space is employed. To construct approximation in time, regularized two-
level schemes are used. The numerical implementation is based on solving
the equation with  the square root of the elliptic operator using an
auxiliary Cauchy problem for a pseudo-parabolic equation.
The scheme of the second-order accuracy in time is based on a 
regularization of the three-level explicit Adams scheme. 
More general problems for the equation with convective terms are considered, too.
The results of numerical experiments are presented for a model two-dimensional problem.
\end{abstract}


\end{frontmatter}

\section{Introduction}

Non-local applied mathematical models based on the use of fractional derivatives in time and space 
are actively discussed in the literature \cite{baleanu2012fractional,kilbas2006theory}. 
Many models, which are used in applied physics, biology, hydrology and finance,
involve both sub-diffusion (fractional in time) and super-diffusion (fractional in space) operators. 
Super-diffusion problems are treated as evolutionary problems with a fractional power of an elliptic operator.
For example, suppose that in a bounded domain $\Omega$ on the set of functions 
$u(\bm x) = 0, \ \bm x \in \partial \Omega$, 
there is defined the operator $\mathcal{D}$: $\mathcal{D}  u = - \triangle u, \ \bm x \in \Omega$. 
We seek the solution of the Cauchy problem for the equation with the square root of an elliptic operator:
\[
 \frac{d u}{d t} + \mathcal{D}^{1/2} u = f(t),
 \quad 0 < t \leq T, 
\] 
\[
 u(0) = u_0,
\] 
for a given $f(\bm x, t)$, $u_0(\bm x), \ \bm x \in \Omega$ using the notation $f(t) = f(\cdot,t)$. 

To solve problems with the square root of an elliptic operator, we can apply
finite volume and finite element methods oriented to using arbitrary
domains and irregular computational grids \cite{KnabnerAngermann2003,QuarteroniValli1994}.
The computational realization is associated with the implementation of the matrix function-vector multiplication.
For such problems, different approaches \cite{higham2008functions} are available.
Problems of using Krylov subspace methods with the Lanczos approximation
when solving systems of linear equations associated with
the fractional elliptic equations are discussed, e.g., in \cite{ilic2009numerical}.
A comparative analysis of the contour integral method, the extended Krylov subspace method, and the preassigned 
poles and interpolation nodes method for solving
space-fractional reaction-diffusion equations is presented in \cite{burrage2012efficient}.
The simplest variant is associated with the explicit construction of the solution using the known
eigenvalues and eigenfunctions of the elliptic operator with diagonalization of the corresponding matrix
\cite{bueno2012fourier,ilic2005numerical,ilic2006numerical}. 
Unfortunately, all these approaches demonstrates too high computational complexity for multidimensional problems.

We have proposed \cite{vabishchevich2014numerical} a computational algorithm for solving
an equation with fractional powers of elliptic operators on the basis of
a transition to a pseudo-parabolic equation.
For the auxiliary Cauchy problem, the standard two-level schemes are applied.
The computational algorithm is simple for practical use, robust, and applicable to solving
a wide class of problems. A small number of pseudo-time steps is required to reach a steady solution.
This computational algorithm for solving equations with fractional powers of operators
is promising when considering transient problems. 

To solve numerically evolutionary equations of first order, as a rule, two-level difference schemes are used 
for approximation in time.
Investigation of stability for such schemes in the corresponding finite-dimensional
(after discretization in space) spaces is based on the general theory of operator-difference schemes
\cite{Samarskii1989,SamarskiiMatusVabischevich2002}. In particular, the backward Euler scheme and Crank-Nicolson 
scheme are unconditionally stable for a non-negative operator. 
As for one-dimensional problems for the space-fractional diffusion equation,
an analysis of stability and convergence for this equation was conducted in \cite{jin2014error} using
finite element approximation in space. A similar study for the Crank-Nicolson scheme was conducted earlier
in \cite{tadjeran2006second} using finite difference approximations in space.  
We separately highlight the works \cite{huang2008finite,sousa2012second,meerschaert2004finite}, 
where the numerical methods for solving one-dimensional in spaces transient problems for a convection and 
space-fractional diffusion equation are considered. 

In the present paper, we construct unconditionally stable two- and three-level schemes for 
the approximate solution of unsteady problems with the square root 
of an elliptic operator. A transition to a new temporary level involves 
solving the standard elliptic problem for the desired solution. 
The main computational costs are associated with the evaluation of the right-hand side 
containing the square root of an elliptic operator.

The paper is organized as follows.
The formulation of an unsteady problem containing the square root
of an elliptic operator is given in Section 2. 
Finite-element approximation in space is discussed in Section 3. 
In Section 4, we construct a regularized difference scheme 
and investigate its stability. Regularized schemes of the second-order accuracy are developed in section 5. 
More general problems of convection-diffusion are discussed in section 6.
The results of numerical experiments are described in Section 7.
The computational algorithm for solving the equation with a fractional 
power of an operator based on the Cauchy problem for a pseudo-parabolic equation is proposed.
At the end of the work the main results of our study are summarized.

\section{Problem formulation}

In a bounded polygonal domain $\Omega \subset R^m$, $m=1,2,3$ with the Lipschitz continuous boundary $\partial\Omega$,
we search the solution for the problem with a fractional power of an elliptic operator.
Define the elliptic operator as
\begin{equation}\label{1}
  \mathcal{D}  u = - {\rm div}  k({\bm x}) {\rm grad} \, u + c({\bm x}) u
\end{equation} 
with coefficients $0 < k_1 \leq k({\bm x}) \leq k_2$, $c({\bm x}) \geq 0$.
The operator $\mathcal{D}$ is defined on the set of functions $u({\bm x})$ that satisfy
on the boundary $\partial\Omega$ the following conditions:
\begin{equation}\label{2}
  k({\bm x}) \frac{\partial u }{\partial n } + \mu ({\bm x}) u = 0,
  \quad {\bm x} \in \partial \Omega ,
\end{equation} 
where $\mu ({\bm x}) \geq \mu_1 > 0, \  {\bm x} \in \partial \Omega$.

In the Hilbert space $H = L_2(\Omega)$, we define the scalar product and norm in the standard way:
\[
  <u,v> = \int_{\Omega} u({\bm x}) v({\bm x}) d{\bm x},
  \quad \|u\| = <u,u>^{1/2} .
\] 
For the spectral problem
\[
 \mathcal{D}  \varphi_k = \lambda_k \varphi_k, 
 \quad \bm x \in \Omega , 
\] 
\[
  k({\bm x}) \frac{\partial  \varphi_k}{\partial n } + \mu ({\bm x}) \varphi_k = 0,
  \quad {\bm x} \in \partial \Omega , 
\] 
we have 
\[
 \lambda_1 \leq \lambda_2 \leq ... ,
\] 
and the eigenfunctions  $ \varphi_k, \ \|\varphi_k\| = 1, \ k = 1,2, ...  $ form a basis in $L_2(\Omega)$. Therefore, 
\[
 u = \sum_{k=1}^{\infty} (u,\varphi_k) \varphi_k .
\] 
Let the operator $\mathcal{A}$ be defined in the following domain:
\[
 D(\mathcal{D} ) = \{ u \ | \ u(\bm x) \in L_2(\Omega), \ \sum_{k=0}^{\infty} | (u,\varphi_k) |^2 \lambda_k < \infty \} .
\] 
Under these conditions the operator $\mathcal{D}$ is self-adjoint and positive defined: 
\begin{equation}\label{3}
  \mathcal{D}  = \mathcal{D} ^* \geq \delta I ,
  \quad \delta > 0 ,    
\end{equation} 
where $I$ is the identity operator in $H$. For $\delta$, we have $\delta = \lambda_1$.
In applications, the value of $\lambda_1$ is unknown (the spectral problem must be solved).
Therefore, we assume that $\delta \leq \lambda_1$ in (\ref{3}).
Let us assume for the square root of the  operator $\mathcal{D}$
\[
 \mathcal{D}^{1/2} u =  \sum_{k=0}^{\infty} (u,\varphi_k) \lambda_k^{1/2}  \varphi_k .
\] 
We seek the solution of the Cauchy problem for the evolutionary 
first-order equation with the square root of the operator $\mathcal{D}$. 
The solution $u(\bm x,t)$ satisfies the equation
\begin{equation}\label{4}
  \frac{d u}{d t} + \mathcal{D}^{1/2} u = f(t),
  \quad 0 < t \leq T,  
\end{equation} 
and the initial condition
\begin{equation}\label{5}
  u(0) = u_0 .
\end{equation} 
The key issue in the study of the computational algorithm for solving the Cauchy problem (\ref{4}), (\ref{5}) 
is to establish the stability of the approximate solution with respect to small perturbations of the initial 
data and the right-hand side in various norms.

\section{Discretization in space}\label{sec:3}

To solve numerically the problem (\ref{4}), (\ref{5}), we employ finite-element 
approximations in space \cite{brenner2008mathematical,Thomee2006}. 
For (\ref{1}) and (\ref{2}), we define the bilinear form
\[
 d(u,v) = \int_{\Omega } \left ( k \, {\rm grad} \, u \, {\rm grad} \, v + c \, u v \right )  d {\bm x} +
 \int_{\partial \Omega } \mu \, u v d {\bm x} .
\] 
By (\ref{3}), we have
\[
d(u,u) \geq \delta \|u\|^2 .  
\]
Define a subspace of finite elements $V^h \subset H^1(\Omega)$.
Let $\bm x_i, \ i = 1,2, ..., M_h$ be triangulation points for the domain $\Omega$.
Define pyramid function $\chi_i(\bm x) \subset V^h, \ i = 1,2, ..., M_h$, where
\[
 \chi_i(\bm x_j) = \left \{
 \begin{array}{ll} 
 1, & \mathrm{if~}  i = j, \\
 0, & \mathrm{if~}  i \neq  j .
 \end{array}
 \right . 
\] 
For $v \in V_h$, we have
\[
 v(\bm x) = \sum_{i=i}^{M_h} v_i \chi_i(\bm x),
\] 
where $v_i = v(\bm x_i), \ i = 1,2, ..., M_h$.

We define the discrete elliptic operator $D$ as
\[
<D y, v> = d(y,v),
\quad \forall \ y,v \in V^h . 
\]
The square root of the operator $D$ is defined similarly to $\mathcal{D}^{1/2}$.
For the spectral problem
\[
D \widetilde{\varphi}_k = \widetilde{\lambda}_k 
\] 
we have
\[
\widetilde{\lambda}_1 \leq \widetilde{\lambda}_2 \leq ... \leq  \widetilde{\lambda}_{M_h},
\quad \| \widetilde{\varphi}_k\| = 1,
\quad k = 1,2, ..., M_h . 
\]
The domain of definition for the operator $D$ is
\[
D(D) = \{ y \ | \ y \in V^h, \ \sum_{k=0}^{M_h} | (y,\widetilde{\varphi}_k) |^2 \widetilde{\lambda}_k < \infty \} .
\]
The operator $D$ acts on a finite dimensional space $V^h$ defined on  the domain $D(D)$
and, similarly to (\ref{3}), 
\begin{equation}\label{6}
D = D^* \geq \delta I ,
\quad \delta > 0 , 
\end{equation} 
where $\delta \leq \lambda_1 \leq \widetilde{\lambda}_1$. 
For the square root of the operator $D$, we suppose
\[
D^{1/2} y = \sum_{k=1}^{M_h} (y, \widetilde{\varphi}_k) \widetilde{\lambda}_k^{1/2} 
\widetilde{\varphi}_k .
\] 
For the problem (\ref{4}), (\ref{5}), we put into the correspondence the operator equation 
for $w(t) \in V^h$:
\begin{equation}\label{7}
 \frac{d w}{d t} + D^{1/2} w = \psi(t), 
 \quad 0 < t \leq T, 
\end{equation} 
\begin{equation}\label{8}
 w(0) = w_0, 
\end{equation} 
where $\psi(t) = P f(t)$, $w_0 = P u_0$ with $P$ denoting $L_2$-projection onto $V^h$.

Now we obtain an elementary a priori estimate for the solution of (\ref{2}), (\ref{3})
assuming that the solution of the problem, coefficients of the elliptic operator, 
the right-hand side and initial conditions are sufficiently smooth.

Let us multiply equation (\ref{2}) by $w$ and integrate it over the domain $\Omega$:
\[
 \left <\frac{d w}{d t}, w \right > + < D^{1/2} w, w > \, = \,
 < \psi, w > .
\]
In view of the self-adjointness and positive definiteness of the operator $D^{1/2}$,  
the right-hand side can be evaluated by the inequality
\[
 < \psi, w > \, \leq \, <D^{1/2} w, w > + \frac{1}{4}  <D^{-{1/2}} \psi, \psi  > .
\] 
By virtue of this, we have
\[
 \frac{d}{d t} \|w\|^2 \leq \frac{1}{2} \|\psi \|^2_{D^{-{1/2}}} ,
\] 
where $\|\psi \|_{D^{-{1/2}}} = <D^{-{1/2}} \psi, \psi >^{1/2}$. 
The latter inequality leads us to the desired a priori estimate:
\begin{equation}\label{9}
 \|w(t)\|^2 \leq \|w_0\|^2 + \frac{1}{2} \int_{0}^{t}\|\psi(\theta) \|^2_{D^{-{1/2}}} d \theta .
\end{equation} 
We will focus on the estimate (\ref{9})
for the stability of the solution with respect to the initial data and the right-hand side 
in constructing discrete analogs of the problem (\ref{7}), (\ref{8}).

\section{Two-level scheme}

Let $\tau$ be a step of a uniform grid in time such that $w^n = w(t^n), \ t^n = n \tau$,
$n = 0,1, ..., N, \ N\tau = T$.
To solve numerically the problem (\ref{7}), (\ref{8}), we use the simplest implicit two-level scheme. 
We start  with the simplest explicit scheme
\begin{equation}\label{10}
 \frac{w^{n+1} - w^{n}}{\tau } + D^{1/2} w^{n} = \psi^{n},
 \quad n = 0,1, ..., N-1,
\end{equation} 
\begin{equation}\label{11}
 w^0 = w_0 .
\end{equation}
When considering the standard evolutionary problems, advantages and disadvantages of explicit schemes 
are well-known (see, e.g., \cite{Samarskii1989,SamarskiiMatusVabischevich2002}).
For problems with the square root of the operator, the main drawback (conditional stability) 
for explicit schemes remains. In addition, the advantage in terms of implementation simplicity is no more.
In this case, the approximate solution at a new time level is determined via (\ref{10}) as
\begin{equation}\label{12}
 w^{n+1} = w^{n} - \tau D^{1/2} w^{n} + \tau \psi^{n} . 
\end{equation} 
Thus, we must calculate $D^{1/2} w^{n}$.

Consider now the possibility of using implicit schemes.
Let us approximate equation (\ref{7}) by the backward Euler scheme:
\begin{equation}\label{13}
 \frac{w^{n+1} - w^{n}}{\tau } + D^{1/2} w^{n+1} = \psi^{n+1},
 \quad n = 0,1, ..., N-1.
\end{equation}
The main advantage of the implicit scheme (\ref{13}) in comparison with (\ref{10})
is its absolute stability. Let us derive for this scheme  the corresponding estimate for stability.

Multiplying equation (\ref{13}) scalarly by $\tau w^{n+1}$, we obtain
\begin{equation}\label{14}
\begin{split}
 < w^{n+1}, w^{n+1}> \ & + \ \tau < D^{1/2} w^{n+1}, w^{n+1}> \  \\
 & = \ < w^{n}, w^{n+1}> + \tau < \psi^{n+1}, w^{n+1}> .
\end{split} 
\end{equation}
The terms on the right side of (\ref{14}) are estimated using the inequalities:
\[
 < w^{n}, w^{n+1}> \ \leq \frac{1}{2} < w^{n+1}, w^{n+1}> + \frac{1}{2} < w^{n}, w^{n}> ,
\] 
\[
 < \psi^{n+1}, w^{n+1}> \ \leq \ < D^{1/2} w^{n+1}, w^{n+1}> + \frac{1}{4}  < D^{-{1/2}} \psi^{n+1}, \psi^{n+1}> .
\]
The substitution into (\ref{14}) leads to the following level-wise estimate:
\[
 \|w^{n+1}\|^2 \leq \|w^{n}\|^2 + \frac{1}{2} \tau \|\psi^{n+1}\|^2_{D^{-{1/2}}} .
\]
This implies the desired estimate for stability:
\begin{equation}\label{15}
 \|w^{n+1}\|^2 \leq \|w_0\|^2 + \frac{1}{2} \sum_{k=0}^{n}\tau \|\psi^{k+1}\|^2_{D^{-{1/2}}} ,
\end{equation} 
which is a discrete analog of the estimate (\ref{9}). 
To obtain the solution at the new time level, it is necessary to solve the problem
\[
 (I +  \tau D^{1/2}) w^{n+1} = w^{n} + \tau \psi^{n} .
\]
In our case, we must calculate the values of $\varPhi(D) b$ for $\varPhi(z) = (1+ \tau z^{1/2})^{-1}$. 

A more complicated situation arises in the implementation of the Crank-Nicolson scheme:
\[
 \frac{w^{n+1} - w^{n}}{\tau } + D^{1/2}  \frac{w^{n+1} + w^{n}}{2} = \psi^{n+1/2},
 \quad n = 0,1, ..., N-1.
\]
In this case, we have 
\[
 \left (I +  \frac{\tau}{2}  D^{1/2} \right ) w^{n+1} = w^{n} - \frac{\tau}{2}  D^{1/2} w^{n}
 + \tau \psi^{n+1/2} ,
\]
i.e., we need to evaluate both $\varPhi(z) = (1+ 0.5 \tau z^{1/2})^{-1}$ and $\varPhi(z) = z^{1/2}$.

The numerical implementation of the above-mentioned approximations in time for the standard parabolic problems in (\ref{7}) is based on calculating the values of $\varPhi(D) b$ 
for $\varPhi(z) = (1+ \sigma \tau z)^{-1}, \ \sigma =0.5, 1$ and $\varPhi(z) = z$. 
For problems with the square root of elliptic operators, 
we apply the approach proposed early in our paper \cite{vabishchevich2014numerical}.
It is based on the computation of $\varPhi(D) b$ for $\varPhi(z) = z^{-{1/2}}$.

For the explicit approximation in time, we rewrite (\ref{12}) in the form
\[
 w^{n+1} = w^{n} - \tau D D^{-{1/2}} w^{n} + \tau \psi^{n} .
\]
Therefore, the computational implementation is based on the evaluation of 
$\varPhi(D) b$ for $\varPhi(z) = z^{-{1/2}}$ and $\varPhi(z) = z$. 
A similar approach is not valid for the  backward Euler scheme (\ref{11}), (\ref{13}) 
and moreover for the Crank-Nicolson scheme. To construct
a more appropriate from computational point of view 
approximations in time for the Cauchy problem (\ref{7}), (\ref{8}), we apply 
the principle of regularization for operator-difference schemes proposed by A.A. Samarskii \cite{Samarskii1989}.

For a regularizing operator $R = R^* > 0$, the simplest regularized scheme 
for solving (\ref{7}), (\ref{8}) has the form (see, e.g., \cite{Vabishchevich2014}):
\begin{equation}\label{16}
 (I + R) \frac{w^{n+1} - w^{n}}{\tau } + D^{1/2} w^{n} = \psi^{n+1},
 \quad n = 0,1, ..., N-1.
\end{equation}
Now we will derive the stability conditions for the regularized scheme (\ref{11}), (\ref{16})  
and after that we will select the appropriate regularizing operator $R$ itself.

Rewrite equation (\ref{16}) in the form
\[
 \left (I +  R - \frac{\tau}{2}D^{1/2} \right ) \frac{w^{n+1} - w^{n}}{\tau } 
 + D^{1/2} \frac{w^{n+1} + w^{n}}{2} = \psi^{n+1} .
\] 
Multiplying it scalarly by $\tau (w^{n+1} + w^{n})$, we get
\[
\begin{split}
 < G w^{n+1}, w^{n+1}> - < G w^{n}, w^{n}> & +  \frac{\tau }{2} <D^{1/2} (w^{n+1} + w^{n}), w^{n+1} + w^{n} > \\
 & = \tau <\psi^{n+1}, w^{n+1} + w^{n} > ,
\end{split}
\] 
where
\begin{equation}\label{17}
 G = I + R - \frac{\tau}{2}D^{1/2} .
\end{equation} 
For $G = G^* > 0$, and $G = I + \mathit{O}(\tau)$,
we obtain the inequality
\[
 \|w^{n+1}\|_G^2 \leq \|w^{n}\|_G^2 + \frac{1}{2} \tau \|\psi^{n+1}\|^2_{D^{-{1/2}}} .
\]
Thus, for the regularized difference scheme (\ref{11}), (\ref{16})
the following estimate for stability with respect to the initial data 
and the right-hand side holds:
\begin{equation}\label{18}
 \|w^{n+1}\|_G^2 \leq \|w_0\|_G^2 + \frac{1}{2} \sum_{k=0}^{n}\tau \|\psi^{k+1}\|^2_{D^{-{1/2}}} .
\end{equation} 

To select an appropriate regularizing operator $R$, we should take into account two conditions, i.e.,
first, to satisfy the inequality $G > 0$, and secondly, to simplify calculations. 
In the scheme (\ref{16}), we put
\begin{equation}\label{19}
 R = \sigma \tau (D + I) .
\end{equation}
In view of
\[
 D - 2 D^{1/2} + I \geq 0,
\] 
in selecting (\ref{19}), we have
\[
 G = I + \sigma \tau (D + I) - \frac{\tau}{2}D^{1/2} \geq 
 I + \left (\sigma - \frac{1}{4} \right ) \tau (D + I) ,
\] 
for $\sigma \geq 0.25$ the inequality  $G > 0$ holds.
The result of our analysis is the following statement.

\begin{thm}\label{Th1}
The regularized scheme (\ref{11}), (\ref{16})  with the regularizer $R$
selected according to (\ref{19}) is unconditionally stable for $\sigma \geq 0.25$.
The approximate solution satisfies the a priori estimate (\ref{17}), (\ref{18}). 
\end{thm}

The transition to a new time level is performed via the formula
\[
 ((1+\sigma \tau) I +  \sigma \tau D ) w^{n+1} =  
 ((1+\sigma \tau) I +  \sigma \tau D) w^{n}  
 - \tau D D^{-{1/2}} w^{n} + \tau \psi^{n+1} . 
\]
Therefore, it is necessary to calculate the values 
$\varPhi(D) b$ for $\varPhi(z) = (1+\tau\widetilde{\sigma} z)^{-1}$ 
and $\varPhi(z) = z^{-{1/1}}$.

\section{The scheme of the second-order accuracy} 

The regularized scheme (\ref{11}), (\ref{16}) has the first-order accuracy in time. Let us consider 
the possibility to construct second-order regularized schemes in the class of three-level schemes.

Among multilevel difference schemes we can select the explicit
Adams methods. To start with the generating scheme, we take the explicit three-level scheme
\begin{equation}\label{20}
 \frac{w^{n+1} - w^{n}}{\tau } + D^{1/2} \frac {3w^n - w^{n-1}}{2} = \psi^{n+1/2},
 \quad n = 1,2, ..., N-1,
\end{equation} 
with specified $w^0, w^1$. It is possible to define $w^1$ via the explicit scheme
\[
 \frac{w^{1} - w^{0}}{\tau } + \left (D^{1/2} - \frac{\tau}{2} D \right ) w^0 =
 \varphi^0 .
\] 
First of all, we formulate a condition for the stability of the scheme.

The scheme (\ref{20}) may be written in the canonical form
for three-level difference schemes \cite{Samarskii1989,SamarskiiMatusVabischevich2002}: 
\begin{equation}\label{21}
  \tilde{B} \frac {w^{n+1} - w^{n-1}}{2\tau} +
  \tilde{R} (w^{n+1} - 2 w^n + w^{n-1}) +
  \tilde{A} w^n = \varphi^n .
\end{equation}
Taking into account
\[
  w^{n+1} - w^{n} =
  \frac{1}{2} (w^{n+1} - w^{n-1}) +
  \frac{1}{2} (w^{n+1} - 2 w^{n} + w^{n-1}),
\]
\[
  3w^{n} - w^{n-1} =
  \frac{1}{2} (w^{n+1} - w^{n-1}) -
  \frac{1}{2} (w^{n+1} - 2 w^{n} + w^{n-1}) + 2w^n,
\]
the scheme (\ref{20}) obtains the representation (\ref{21}), where
\[
   \tilde{B} = I + \frac \tau 2 D^{1/2},
   \quad \tilde{R} = \frac 1{2\tau} \left (I - \frac \tau 2 D^{1/2} \right ) ,
   \quad \tilde{A} = D^{1/2}.
\]
The necessary and sufficient condition for stability of the three-level scheme
(\ref{20}) with self-adjoint operators (see \cite{Samarskii1989,SamarskiiMatusVabischevich2002}) has the form
\begin{equation}\label{22}
   \tilde{B} \ge 0,
   \quad \tilde{R} \ge \frac 14 \tilde{A},
   \quad \tilde{A} > 0 .
\end{equation}
By (\ref{22}), the explicit scheme (\ref{20}) is stable under the condition
\[
  \tilde{R} - \frac{1}{4} \tilde{A} =
  \frac{1}{2\tau} I - \frac{1}{2} D^{1/2} \ge 0.
\]
This gives the constraint on a time step:
\[
  \tau \le \frac{1}{\|D^{1/2}\|} .
\]
In compare with the two-level explicit scheme, this reduce the permissible time step by a factor of two.

The construction of a three-level regularized scheme is conducted similarly to the scheme (\ref{11}), (\ref{16}).
Remaining in the class of the schemes with the second-order approximation, instead of (\ref{22}), we use
the regularized scheme
\begin{equation}\label{23}
 (I + S) \frac{w^{n+1} - w^{n}}{\tau } + D^{1/2} \frac {3w^n - w^{n-1}}{2} = \psi^{n+1.2},
 \quad n = 1,2, ..., N-1.
\end{equation}

\begin{thm}\label{Th2}
The regularized scheme (\ref{23}) with the regularizer $S = \sigma \tau^2 D$ is unconditionally stable 
for $\sigma \geq 1/4$.
\end{thm}

\begin{pf}
The scheme (\ref{23}) may be written in the canonical form (\ref{21}) with
\[
   \tilde{B} = I + S + \frac \tau 2 D^{1/2},
   \quad \tilde{R} = \frac 1{2\tau} \left (I + S - \frac \tau 2 D^{1/2} \right ) ,
   \quad \tilde{A} = D^{1/2}.
\]
For $S \geq 0$ the condition for stability of (\ref{22}) holds if
\[
\begin{split}
  \frac{1}{2 \tau} \left (I + S - \frac{\tau }{2} D^{1/2} \right )  - \frac 14 D^{1/2} & = \frac{1}{2 \tau} \left (I + \sigma \tau^2 D - \frac{\tau }{2} D^{1/2} \right )  - \frac 14 D^{1/2}  \\
  & = \frac{1}{2 \tau} \left ( \left (I - \frac{\tau }{2} D^{1/2} \right )^2 +
  \left (\sigma - \frac 14 \right ) \tau^2 D \right ) .  
\end{split}
\]
Thus, the regularized scheme (\ref{23}) is unconditionally
stable if $\sigma \ge 1/4$ in (\ref{19}) .
\end{pf}

The computational implementation of the scheme (\ref{19}), (\ref{23}) 
is similar to the numerics of a two-level scheme (\ref{16}), (\ref{19}).

\section{Convection-diffusion problems} 

Convection-diffusion problems are typical for mathematical models of fluid mechanics.
Heat transfer as well as impurities spreading are occurred not only due to diffusion, 
but result also from medium motion. Principal features of physical and chemical phenomena
observing in fluids and gases \cite{landau_fluid_1987,batchelor_introduction_2000}
are generated by media motion resulting from various forces.
Computational algorithms for the numerical solution of such problems are of great importance;
they are discussed in many publications (see, e.g., \cite{HundsdorferVerwer2003,MortonKellogg1996}).
Here we  discuss the problem of constructing unconditionally stable schemes for 
convection-diffusion problems, where diffusive transport is described by the square root 
of an elliptic operator.

The convection is provided by transient field of medium velocity  $\bm v (\bm x, t)$.
On the boundary the following condition is satisfied
\begin{equation}\label{24}
 \bm v (\bm x, t) = 0, \quad \bm x \in \partial \Omega . 
\end{equation}  
We take the convection operator in the skew-symmetric form \cite{Vabishchevich2014,ChurbanovVabishchevich2013}:
\begin{equation}\label{25}
  \mathcal{C}  u = {\rm div}  (\bm v \, u) - \frac{1}{2} {\rm div} \bm v  \, u .
\end{equation} 
Using the restrictions  (\ref{24}) the convection operator is skew-symmetric:
\begin{equation}\label{26}
  \mathcal{C}  = -  \mathcal{C}^* .
\end{equation} 

We seek the solution of the Cauchy problem for the convection-diffusion equation.
The solution $u(\bm x,t)$ satisfies the equation
\begin{equation}\label{27}
  \frac{d u}{d t} +  \mathcal{C} u + \mathcal{D}^{1/2} u = f(t),
  \quad 0 < t \leq T,  
\end{equation} 
and the initial condition (\ref{5}).
The main feature of the problem is that the evolutionary 
first-order equation has the square root of the operator $\mathcal{D}$. 

The discrete convection operator is introduced similarly.
Taking into account (\ref{24}), (\ref{25}) we  set
\[
 c(y,w) = \int_{\Omega } {\rm div} (\bm v \, y ) \, w \,  d {\bm x} 
 - \frac{1}{2} \int_{\Omega } {\rm div} \bm v \, y w \, d {\bm x} .
\] 
We define the operator $C$ as
\[
 c(y,w) = \ <Cy, w>,
 \quad \forall \ y,w \subset V^h .
\] 
Similarly to (\ref{26}) we have
\begin{equation}\label{28}
 C = - C^* .
\end{equation} 

For the problem (\ref{5}), (\ref{27}), we put into the correspondence the operator equation 
for $w(t) \subset V^h$:
\begin{equation}\label{29}
 \frac{d w}{d t} + C w + D^\alpha w = \psi(t), 
 \quad 0 < t \leq T, 
\end{equation}
taking into account the previously introduced notation and the initial condition (\ref{8}).

Let us represent the two-level regularized scheme (\ref{8}), (\ref{29}) in the form
\begin{equation}\label{30}
 (I + R) \frac{w^{n+1} - w^{n}}{\tau } + C \frac{w^{n+1} + w^{n}}{2} +  D^{1/2} w^{n} = \psi^{n+1},
 \quad n = 0,1, ..., N-1.
\end{equation}
The stability of this scheme is considered in the same way as for the scheme (\ref{16}). 
The key moment is the use of the skew-symmetry property of convective transport operator
(\ref{28}). The scheme (\ref{30}) may be written as
\[
 \left (I +  R - \frac{\tau}{2}D^{1/2} \right ) \frac{w^{n+1} - w^{n}}{\tau } 
 + (C +  D^{1/2})\frac{w^{n+1} + w^{n}}{2} = \psi^{n+1} .
\] 
Multiplying it scalarly by $\tau (w^{n+1} + w^{n})$, in view of (\ref{28}), we again obtain
\[
\begin{split}
 < G w^{n+1}, w^{n+1}> - < G w^{n}, w^{n}> & +  \frac{\tau }{2} <D^{1/2} (w^{n+1} + w^{n}), w^{n+1} + w^{n} > \\
 & = \tau <\psi^{n+1}, w^{n+1} + w^{n} > ,
\end{split}
\] 
where operator $G$ is defined according to (\ref{17}). 

\begin{thm}\label{Th3}
The regularized scheme (\ref{11}), (\ref{30})  with the regularizer $R$
selected according to (\ref{19}) is unconditionally stable for $\sigma \geq 0.25$.
The approximate solution satisfies the a priori estimate (\ref{17}), (\ref{18}). 
\end{thm}

At the new time level we solve the problem
\[
 \Big ( (1+\sigma \tau ) I +   \frac{\tau}{2} C +  \sigma \tau D \Big ) w^{n+1} =  
 \Big ( (1+\sigma \tau ) I  -  \frac{\tau}{2} C + \sigma \tau D \Big ) w^{n}  
 - \tau D D^{-{1/2}} w^{n} + \tau \psi^{n+1} . 
\]
Therefore, it is necessary to solve the discrete elliptic problem with the operator
\[
 B = ((1+\sigma \tau ) I +  \frac{\tau}{2} C + \sigma \tau D
\] 
and compute $D^{-\beta} w^{n}$.

\section{Numerical experiments} 

The main peculiarity of solving the  Cauchy problems for equations  (\ref{16}), (\ref{23}) and
(\ref{30}) with the regularizer (\ref{19}) is the necessity to evaluate values
\[
 g^n = D^{-1/2} w^{n}, 
 \quad n = 0,1, ..., N-1 .
\]
The computational algorithm is based on the consideration of the auxiliary 
Cauchy problem \cite{vabishchevich2014numerical}.

Computational implementation of regularized schemes is based on the previously proposed approach involving 
the solution of the auxiliary Cauchy problem for a pseudo-parabolic equation. 
It makes possible to construct an unconditionally stable scheme, to obtain an appropriate estimate of stability, 
and to conduct numerical experiments for the unsteady problem of fractal diffusion. Any other known approaches 
can also be applied.  Here we focus on numerical algorithms with the calculation of $D^{-\alpha} y, \ \alpha = 1/2$. 
But it is possible to employ methods based on an approximation  of the operator $D^{-\alpha}$. 
In this connection, an approach of great interest is to approximate the operator $D^{-\alpha}$ 
using operator terms $(I + \gamma_i D)^{-1}$ \cite{Bonito}.

Assume that
\[
 y(s) = \delta^{{1/2}} (s (D - \delta I) + \delta I)^{-{1/2}} y(0) ,
\]
then for the determination of $g^n$, we can put
\begin{equation}\label{31}
 g^n = y(1),
 \quad y(0) = \delta^{-{1/2}}  w^{n}.  
\end{equation} 
The function $y(s)$ satisfies the evolutionary equation
\begin{equation}\label{32}
  (s G + \delta I) \frac{d y}{d s} + \frac{1}{2}  G y = 0 ,
\end{equation} 
where $G = D - \delta I \geq  0$.
Thus, the calculation of $D^{-{1/2}} w^{n}$  is based on the solution of 
the Cauchy problem (\ref{31}), (\ref{32}) within  the unit interval for the pseudo-parabolic equation.

To solve numerically the problem (\ref{31}), (\ref{32}),
we use the simplest implicit two-level scheme.
Let $\eta$ be a step of a uniform grid in time such that $y_k = y(s_k), \ s_k = k \eta$,
$k = 0,1, ..., K, \ K\eta = 1$.
Let us approximate equation (\ref{32})  by the backward Euler scheme
\begin{equation}\label{33}
  (s_{k+1} G + \delta I) \frac{ y_{k+1} - y_{k}}{\eta } + 
  \frac{1}{2}   G  y_{k+1}  = 0,  
  \quad k = 0,1, ..., K-1,
\end{equation}
\begin{equation}\label{34}
  y_0 = \delta^{-{1/2} } w^{n} . 
\end{equation}
For the Crank-Nicolson scheme, we have
\begin{equation}\label{35}
  (s_{k+1/2} G + \delta I) \frac{ y_{k+1} - y_{k}}{\eta } + 
  \frac{1}{2}   G  \frac{ y_{k+1} + y_{k}}{2}  = 0,  
  \quad k = 0,1, ..., K-1 .
\end{equation}
The difference scheme (\ref{34}), (\ref{35}) approximates the problem 
(\ref{31}), (\ref{32})
with the second order by $\eta $, whereas for scheme (\ref{33}), (\ref{34})  we have only the first order.

The above two-level schemes are unconditionally stable. 
The corresponding level-wise estimate has the form
\begin{equation}\label{36}
 \|y_{k+1}\| \leq \|y_{k}\| ,
 \quad k = 0,1, ..., K-1 . 
\end{equation}
To prove (\ref{36}) (see \cite{vabishchevich2014numerical}), it is sufficient to multiply scalarly
equation (\ref{33}) by $y_{k+1}$ and equation (\ref{35}) by $y_{k+1}+y_{k}$. 
Taking into account (\ref{34}), from (\ref{36}), we obtain
\begin{equation}\label{37}
 \|y_K\| \leq \delta^{-{1/2} } \|w^{n}\| . 
\end{equation} 

\begin{figure}[htp]
  \begin{center}
  \begin{tikzpicture}[scale = 2]
     \path[draw=black, fill=white] (0,0) circle (1.);
     \draw [ultra thick, fill=gray!20] (0,0) -- (0,1) arc [radius=1.0, start angle=90, end angle= 0] -- (1,0) -- (0,0);
     \draw [->] (1,0) -- (1.75,0);
     \draw [->] (0,1) -- (0,1.75);
     \draw(-0.1,-0.1) node {$0$};
     \draw(1.1,-0.1) node {$1$};
     \draw(1.74,-0.1) node {$x_1$};
     \draw(-0.1,1.1) node {$1$};
     \draw(-0.1,1.65) node {$x_2$};
     \draw(0.45,0.45) node {$\Omega$};
     \draw(0.6,-0.13) node {$\Gamma_1$};
     \draw(-0.13,0.6) node {$\Gamma_2$};
     \draw(0.8,0.8) node {$\Gamma_3$};
  \end{tikzpicture}
  \caption{Computational domain $\Omega$}
  \label{f-1}
  \end{center}
\end{figure}
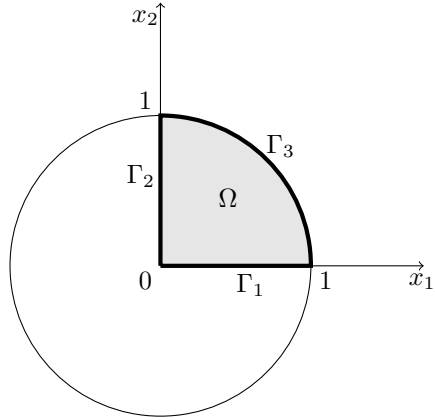

The test problem is constructed using the exact solution of the problem in the unit circle. 
The computational domain is a quarter of the circle (see Fig.~\ref{f-1}).

Let
\[
 \mathcal{D} u = - \Delta u ,
 \quad \bm x \in \Omega , 
\]
under the boundary conditions 
\[
 \frac{\partial u}{\partial n} = 0, 
 \quad \bm x \in \Gamma_1,
 \quad \bm x \in \Gamma_2, 
\] 
\[
 \frac{\partial u}{\partial n} + \mu u = 0,
 \quad \bm x \in \Gamma_1,
 \quad \mu = \const .  
\]
We consider the case, where the solution depends only on 
$r$, and $r = (x_1^2 + x_2^2)^{1/2}$. By virtue of this
\[
 \mathcal{D} u = - \frac{1}{r} \frac{d}{d r} \left ( r  \frac{d u}{d r} \right )  ,
 \quad 0 < r < 1 ,
\] 
\[
 \frac{d u}{d r} + \mu u = 0,
 \quad r = 1 . 
\] 

The solution of the spectral problem
\[
 - \frac{1}{r} \frac{d}{d r} \left ( r  \frac{d \varphi_k}{d r} \right )  = \lambda_k \varphi_k,
 \quad 0 < r < 1 ,
\] 
\[
 \frac{d \varphi_k}{d r} + \mu \varphi_k = 0,
 \quad r = 1 , 
\]
is well-known (see, e.g., \cite{Polyanin,Carslaw}).
Eigenfunctions are represented as zero-order Bessel functions:
\[
 \varphi_k(r) = J_0(\sqrt{\lambda_k} r),
\] 
whereas eigenvalues $\lambda_k = \nu_k^2$, where $ \nu_k, \ k = 1,2, ...$ are roots of the equation
\begin{equation}\label{38}
 \nu  J^{'}_0(\nu ) + \mu J_0(\nu ) = 0 . 
\end{equation}
The general solution of the homogeneous ($f(t) = 0$) Cauchy problem for  equation (\ref{4}) is
\[
 u(r,t) = \sum_{k=1}^{\infty} a_k \exp(-\nu_k t) J_0(\nu_k r) .
\] 

\begin{figure}
  \begin{center}
    \includegraphics[width=0.8\linewidth] {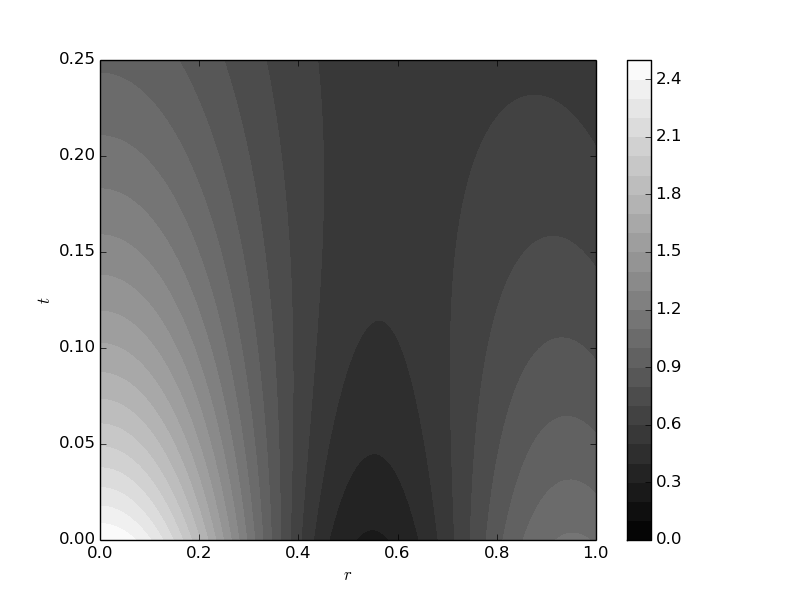}
	\caption{The exact solution for $\mu=1$}
	\label{f-2}
  \end{center}
\end{figure}

\begin{figure}
  \begin{center}
    \includegraphics[width=0.8\linewidth] {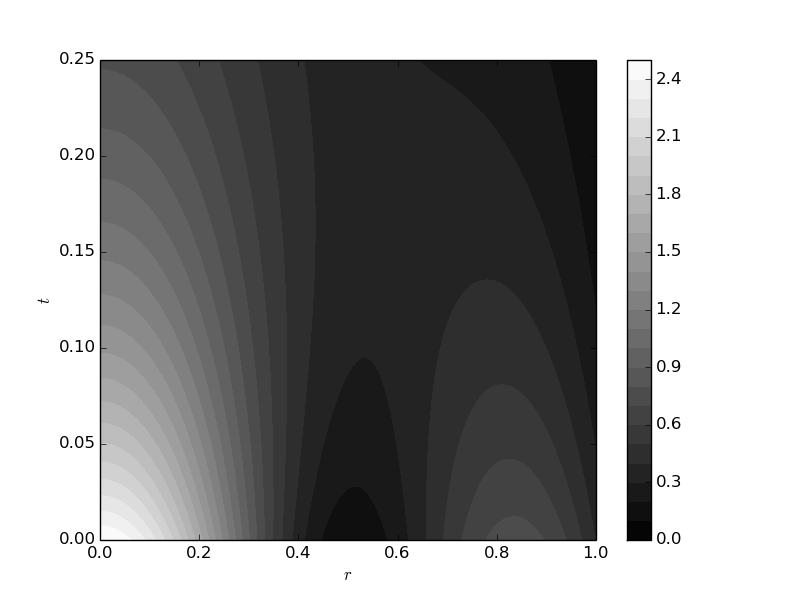}
	\caption{The exact solution for $\mu=10$}
	\label{f-3}
  \end{center}
\end{figure}

\begin{figure}
  \begin{center}
    \includegraphics[width=0.8\linewidth] {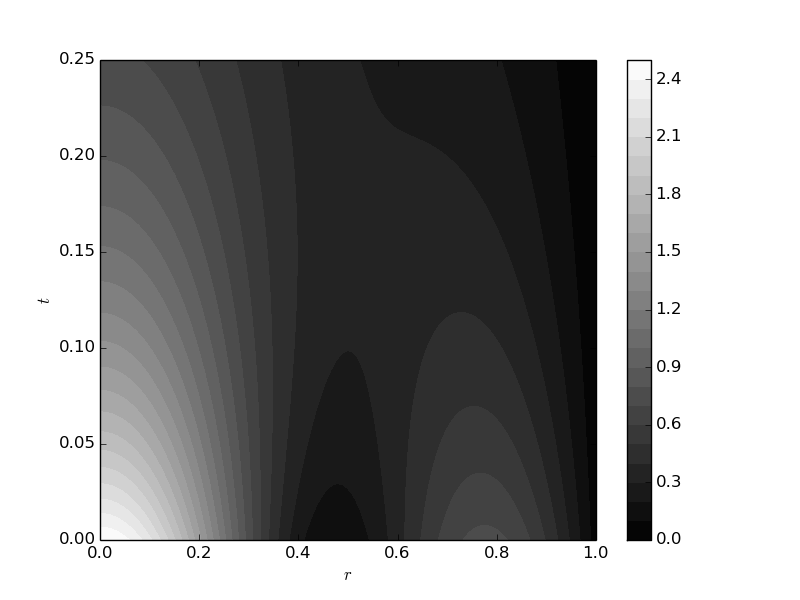}
	\caption{The exact solution for $\mu=100$}
	\label{f-4}
  \end{center}
\end{figure}

To verify the accuracy of the approximate solution of the time-dependent equation with the square root 
of an elliptic operator, we use the exact solution
\begin{equation}\label{39}
 u(r,t) = \exp(-\nu_1 t) J_0(\nu_1 r) + 1.5 \exp(-\nu_3 t) J_3(\nu_3 r),
 \quad r = (x_1^2 + x_2^2)^{1/2} .
\end{equation}
The values of the roots $\nu_1, \ \nu_3$ for different values of the boundary condition $\mu$ are given 
in Table~\ref{tab-1}. The exact solution for $T = 0.25$ at different values of $\mu$ is shown in 
Figs.~\ref{f-2}-\ref{f-4}.

\begin{table}
\begin{center}
 \caption{The roots of equation (\ref{38})}
 \begin{tabular}{cccc}\label{tab-1}
  $k $ &  $\mu = 1$  & $\mu = 10$  & $\mu = 100$  \\
  \hline
  1 & 1.25578371  & 2.17949660  & 2.38090166 \\
  3 & 7.15579917  & 7.95688342  & 8.56783165 \\
\end{tabular}
\end{center} 
\end{table} 

\begin{figure}
  \begin{center}
    \includegraphics[width=0.65\linewidth] {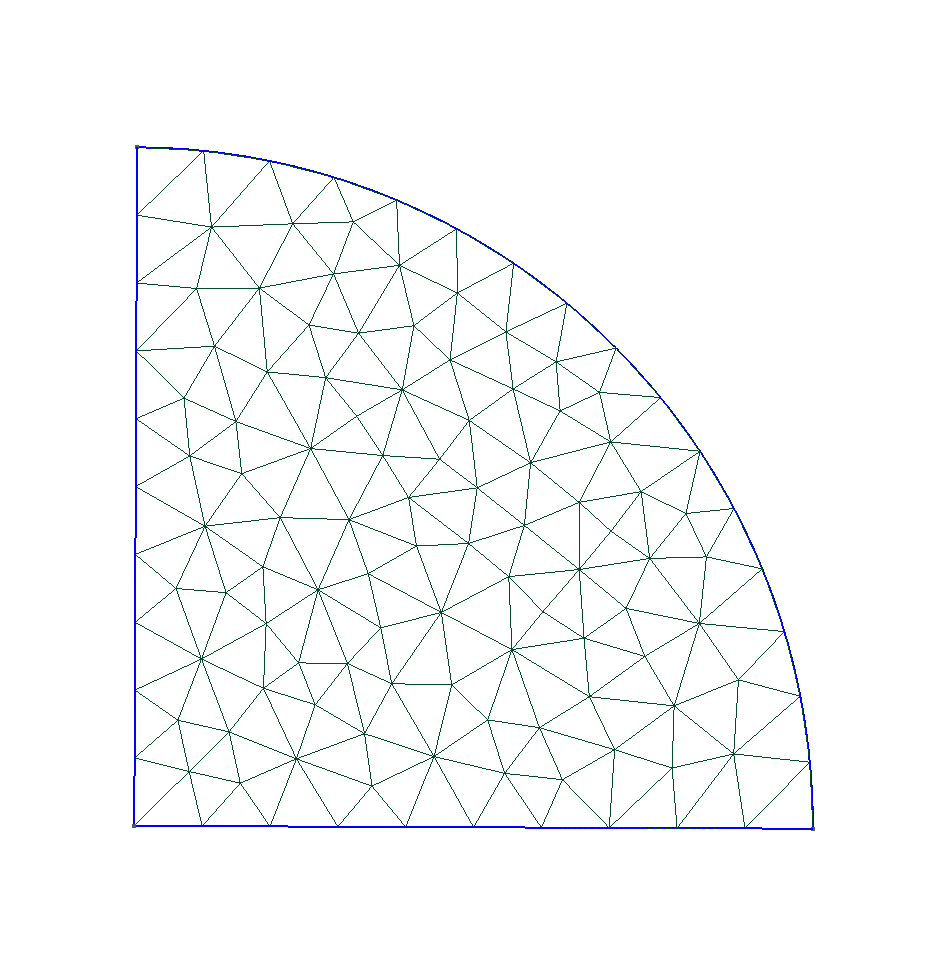}
	\caption{Grid 1: 123 vertices and 208 cells}
	\label{f-5}
  \end{center}
\end{figure}

\begin{figure}
  \begin{center}
    \includegraphics[width=0.65\linewidth] {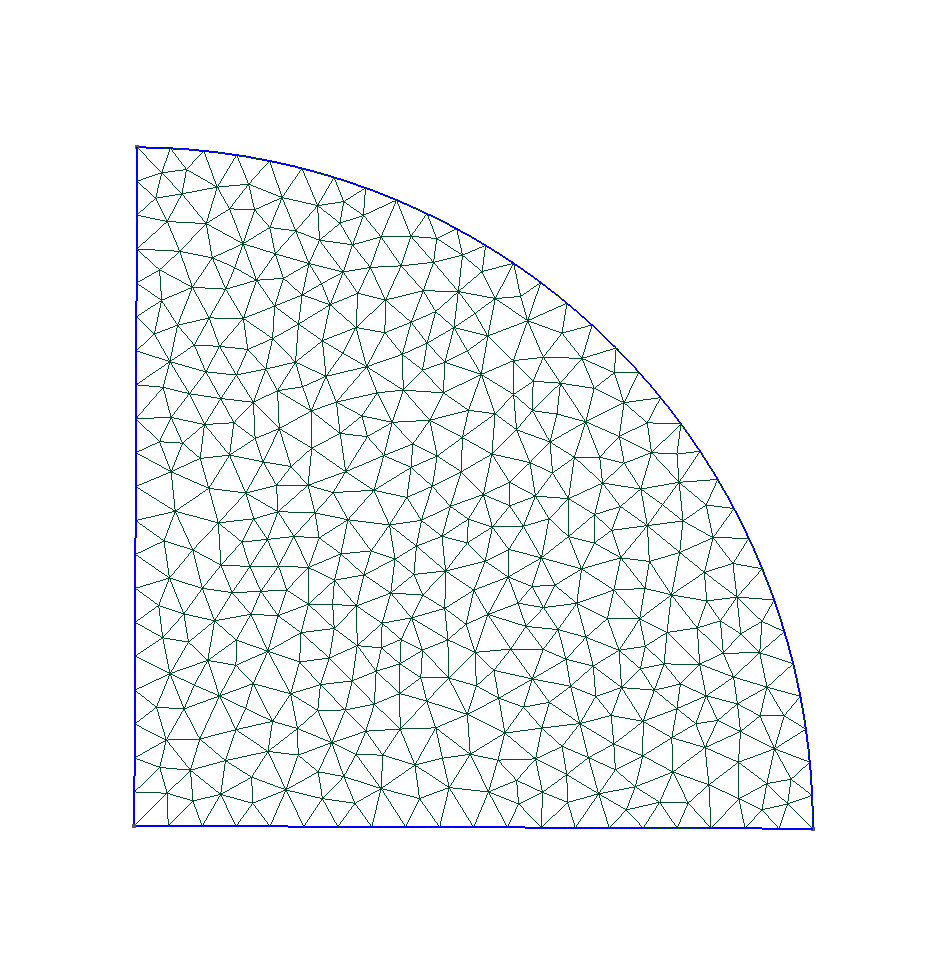}
	\caption{Grid 2: 461 vertices and 848 cells}
	\label{f-6}
  \end{center}
\end{figure}

\begin{figure}
  \begin{center}
    \includegraphics[width=0.65\linewidth] {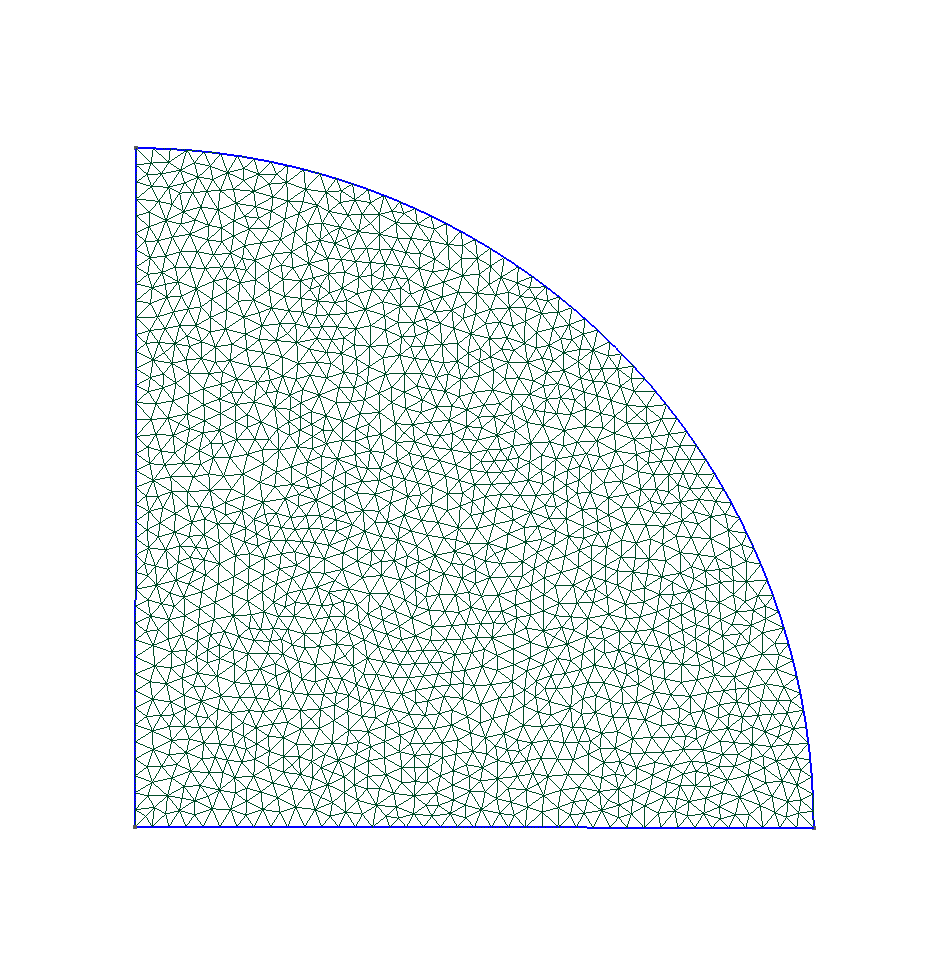}
	\caption{Grid 3: 1731 vertices and 3317 cells}
	\label{f-7}
  \end{center}
\end{figure}

Predictions were performed on a sequence of refining grids, which are shown in Figs.~\ref{f-5}-\ref{f-7}. 
Suppose that the parameter $\delta$ in the estimate (\ref{6}) for $\mu = 1, 10, 100$ is equal to  1. 
The scheme (\ref{35})  is used at each time level on sufficiently fine grid with $K = 100$.
An approximate solution is compared with the exact solution at the final time moment $u(\bm x,T)$. 
Error estimation is performed in $L_2(\Omega)$ and $L_\infty (\Omega)$:
\[
 \varepsilon_2 = \|w^N(\bm x) - u(\bm x,T)\|,
 \quad  \varepsilon_\infty  = \max_{\bm x \in \Omega} |w^N(\bm x) - u(\bm x,T) | .
\] 

We estimated the accuracy of the regularized scheme (\ref{11}), (\ref{16}) on various grids in time and space.
The weighting parameter $\sigma$ was chosen to be  $0.25$.
Table~\ref{tab-2} demonstrates predictions for the model problem with 
$\mu =10$. The dependence of the accuracy on the boundary condition (on the parameter 
$\mu$) is illustrated by Table~\ref{tab-3}. The calculations were performed on grid 2.

\begin{table}[!h]
\caption{Error of the solution  for the regularized scheme (\ref{11}), (\ref{16}) ($\mu =10$)}
\begin{center}
\begin{tabular}{|c|ccccc|} \hline
grid & $N$ & 25 & 50 & 100 & 200  \\ \hline
1    & $\varepsilon_2$       &  0.01316779  &  0.00564075 &  0.00486612  &  0.00616606 \\
     & $\varepsilon_\infty$  &  0.00244349  &  0.00078090 &  0.00114533  &  0.00158964 \\   \hline
2    & $\varepsilon_2$       &  0.01521770  &  0.00784386 &  0.00398968  &  0.00203974 \\
     & $\varepsilon_\infty$  &  0.00375353  &  0.00163994 &  0.00059325  &  0.00020720 \\   \hline
3    & $\varepsilon_2$       &  0.01459601  &  0.00709760 &  0.00332100  &  0.00144008 \\
     & $\varepsilon_\infty$  &  0.00244349  &  0.00199792 &  0.00092803  &  0.00039236 \\   \hline
\end{tabular}
\end{center}
\label{tab-2}
\end{table}

\begin{table}[!h]
\caption{Error of the solution  at various $\mu$}
\begin{center}
\begin{tabular}{|c|ccccc|} \hline
$\mu$   & $N$ & 25 & 50 & 100 & 200  \\ \hline
1    & $\varepsilon_2$       &  0.01192779  &  0.00580199 &  0.00267418  &  0.00157455 \\
     & $\varepsilon_\infty$  &  0.00285443  &  0.00114226 &  0.00034994  &  0.00032851 \\   \hline
10   & $\varepsilon_2$       &  0.01521770  &  0.00784386 &  0.00398968  &  0.00203974 \\
     & $\varepsilon_\infty$  &  0.00375353  &  0.00163994 &  0.00059325  &  0.00020720 \\   \hline
100  & $\varepsilon_2$       &  0.01744919  &  0.00892991 &  0.00447231  &  0.00221753 \\
     & $\varepsilon_\infty$  &  0.00406940  &  0.00179382 &  0.00066321  &  0.00019830 \\   \hline
\end{tabular}
\end{center}
\label{tab-3}
\end{table}

\section{Conclusion} 

\begin{enumerate}
 \item There is considered a nonclassical problem with the initial data, which is described by
 an evolutionary equation of first order with the square root of an elliptic operator. 
 The multidimensional problem is approximated in space using standard 
 finite-element piecewise-linear approximations. A priori estimate of stability with respect to 
 the initial data and the right-hand side is provided.
 \item A regularized two-level scheme of the first-order accuracy in time is constructed.  
 Its unconditional stability is shown. The numerical implementation involves
 the solution of the standard grid elliptic problem at each time level. 
 In addition, in these problems, we need to solve an operator equation 
 for the square root of an elliptic operator.
 \item The possibility of constructing schemes of the second-order accuracy in time 
 for evolutionary problems with the square root of an elliptic operator is studied.
  The improvement of accuracy is achieved via regularization of a three-level scheme.
 \item If convective transport is taken into account, then the skew-symmetric operator is added to the equation. 
 A two-level regularized scheme with the skew-symmetric approximation of convective transport is proposed.
 It demonstrates unconditional stability.
 \item The efficiency of the proposed regularized scheme was demonstrated 
 through the numerical solution of a two-dimensional test problem.
\end{enumerate}

\section*{Acknowledgements}

This work was supported by the Russian Foundation for Basic Research  (projects 14-01-00785, 15-01-00026).

\end{document}